\documentclass[12pt]{amsart}

\usepackage[mathscr]{euscript}
\usepackage{amscd}
\theoremstyle{plain}
\newtheorem{theorem}{Theorem}[section]
\newtheorem{lemma}[theorem]{Lemma}

\theoremstyle{definition}

\numberwithin{equation}{section}

\newcommand{\R}{\mathbb R}
\newcommand{\Sbb}{\mathbb S}
\newcommand{\Hbb}{\mathbb H}

\newcommand{\up}{\Upsilon}
\newcommand{\pa}{\partial}

\newcommand{\gb}{\overline g}
\newcommand{\Xb}{\overline X}
\newcommand{\Yb}{\overline Y}

\newcommand{\rh}{\hat r}
\newcommand{\ep}{\epsilon}
\newcommand{\al}{\alpha}
\newcommand{\be}{\beta}
\newcommand{\ga}{\gamma}
\newcommand{\eh}{\hat{\epsilon}}
\newcommand{\gh}{\hat g}
\newcommand{\om}{\omega}
\newcommand{\Vol}{\mbox{Vol}}
\newcommand{\Area}{\mbox{Area}}
\newcommand{\Ric}{\mbox{Ric}}

\begin{document}

\title[Volume and Area Renormalizations]
{Volume and Area Renormalizations for\\[8pt] Conformally Compact Einstein
Metrics} 

\author[C.R. Graham]{C. Robin Graham}
\address{
  Department of Mathematics\\
  University of Washington \\
  Box 354350
  Seattle, WA 98195-4350}
\email{robin@math.washington.edu}
\maketitle

\thispagestyle{empty}

\section{Introduction}\label{intro}
It has long been known that there are very close connections between the
geometry of hyperbolic space $\Hbb^{n+1}$ of $n+1$ dimensions and the
conformal geometry of the $n$-sphere $\Sbb^n$, viewed as the sphere at
infinity of $\Hbb^{n+1}$.  In recent years it has been realized that it is
fruitful to consider generalizations of some of these connections when
$\Hbb^{n+1}$ is replaced by a ``conformally compact'' Einstein
manifold $X$ of negative scalar curvature, and $\Sbb^n$ is replaced by a
compact conformal manifold $M$, the ``conformal infinity'' of $X$.  Quite
recently there has been a great deal of interest in the physics community
in a correspondence (the so-called Anti-de Sitter/Conformal Field Theory
(AdS/CFT) correspondence) proposed by Maldacena \cite{mal} between string
theory 
and supergravity on such $X$ and supersymmetric 
conformal field theories on $M$.  In this
article we describe some new purely geometric invariants of conformally 
compact Einstein manifolds and of their minimal submanifolds 
which have been discovered via this correspondence. 

The relevant notion of conformal infinity is that introduced by Penrose.
A Riemannian metric $g_+$ on the interior $X^{n+1}$ of a compact manifold
with boundary $\Xb$ is said to be conformally compact if $\gb\equiv r^2g_+$ 
extends continuously (or with some degree of smoothness) as a metric
to $\Xb$, where $r$ is a defining
function for $M=\pa X$, i.e. $r>0$ on $X$ and $r=0$, $dr \neq 0$ on $M$.  
The restriction of $\gb$ to $TM$ rescales upon changing $r$, 
so defines invariantly a conformal class of metrics
on $M$, the conformal infinity of $g_+$.  We are concerned 
with conformally compact metrics $g_+$ which satisfy the Einstein
condition $\Ric(g_+) = -ng_+$.
At least near the hyperbolic metric, these can be
parametrized by their conformal infinities:  in \cite{gl} it is shown that 
each conformal structure on $\Sbb^n$ sufficiently near the standard one is
the conformal infinity of a unique (up to diffeomorphism) conformally
compact Einstein metric on the ball near the hyperbolic metric.

The volume $\Vol(X)$ of any conformally compact manifold $X$ 
is infinite.
An appropriate renormalization of $\Vol(X)$ for $X$ Einstein 
gives rise to the new volume
invariants.  In the physics setting, $\Vol(X)$ arises from 
a concrete procedure outlined by Witten \cite{wit} and independently by
Gubser, Klebanov, and Polyakov \cite{gkp}, following the suggestion of
Maldacena, for calculating
observables in a conformal field theory on $M$ via supergravity and string
theory on $X$.  Under various limits and
approximations, the partition function of a conformal field theory on $M$
is given in terms of the gravitational action on $X$, which for an Einstein
metric $g_+$ is proportional to the volume $\Vol(X)$.

The volume renormalization was carried out by Henningson and Sken\-deris in
\cite{hs1}.  As shown in \cite{fg} and
\cite{gl}, each representative metric on $M$ for the conformal infinity
determines a special defining function $r$ in a neighborhood of $M$.  
As $\ep \rightarrow 0$, the function $\Vol(\{r>\ep\})$ has an asymptotic
expansion in negative powers of $\epsilon$, and a $\log\ep$ term if $n$ is
even.  The coefficients of the negative powers of $\ep$ depend on
the representative conformal metric used to determine $r$.  However, it
turns out that if $n$ is odd, then the constant term in the expansion is
independent of this choice, so is a global invariant of the metric $g_+$.  
If $n$ is even, the constant term is not invariant, giving rise to a
so-called conformal anomaly.  However, in this case the coefficient of the
$\log\ep$ term is invariant, and in fact is given by the integral of a
local curvature expression over $M$.  The $\log\ep$ coefficient is
therefore actually a conformal invariant of $M$, independent of which
$(X,g_+)$ might have been chosen with conformal infinity $M$.

Various of the conformal field theories to which the AdS/CFT
correspondence applies contain observables associated to submanifolds
$N$ of $M$.  According to the correspondence, in a suitable approximation 
the expectation value of such
an observable can be calculated in terms of the area $A(Y)$ in the $g_+$
metric of minimal submanifolds $Y$ of $X$ with $\pa Y =N$. 
Existence theory for such minimal submanifolds is discussed for hyperbolic
$X$ in \cite{a1}, \cite{a2}.  
As in the volume case, necessarily $A(Y) = \infty$, 
so one is led to consideration
of renormalizing the area of a minimal submanifold.  This
renormalization was discussed in hyperbolic space for $\dim N = 1,2$ 
in \cite{bcfm} and in general in \cite{gw}.  If $r$ is the special
defining function associated to a conformal representative on $M$ as above, 
then $\Area(Y\cap \{r>\ep\})$ has an expansion in negative powers of $\ep$,
and again a $\log \ep$ term if $k=\dim(N)$ is even.  The invariance
properties of the coefficients are similar to those above.  If $k$ is odd,
then the constant term in the expansion is independent of the choice of
conformal representative on $M$, so is a global invariant of the minimal
submanifold $Y$.  If $k$ is even, there is a conformal anomaly for the
constant term, but the $\log\ep$ coefficient is a conformal invariant of
the submanifold $N$ of $M$.  One can calculate explicitly the $\log\ep$
coefficient for $k=2$; it turns out to be a version on a general conformal 
manifold of the Willmore functional of a surface in conformally flat
space.  Even in the conformally flat case, this relationship between the
Willmore functional of a surface and the renormalization of the area of
a minimal extension seems to be of some interest.  The Willmore
functional is called the ``rigid string action'' in the physics literature
(\cite{bcfm}, \cite{p}).

In \S2. we review some of the basic properties of conformally compact
Einstein metrics.  
In \S3. we discuss the results of \cite{hs1}: the derivation of the volume
renormalization and resulting invariants and anomalies and the explicit
identification of the $\log\ep$ coefficient and anomaly for $n=4,6$.  We
also calculate the renormalized volume for $\Hbb^{n+1}$ when
$n$ is odd; it turns out that its sign depends on the parity of $(n+1)/2$.
In \S4. we review the area renormalization for minimal submanifolds,
following \cite{gw}.

We remark that in order to justify the derivation of the asymptotic
expansions in $\ep$ of the volume and area, we have to assume that the
Einstein metric $g_+$ and the minimal submanifold $Y$ are sufficiently
regular at infinity.  Here sufficiently regular means that they have 
asymptotic expansions to high enough order, in general involving log terms, 
which formally solve the Einstein or minimal area equations.  One
expects that if the conformal structure on $M$ and the submanifold $N$ are
smooth, then any conformally compact Einstein metric $g_+$ and minimal
submanifold $Y$ will have such regularity, assuming they take on the
boundary data in a suitable sense.   Some regularity results
for minimal submanifolds of hyperbolic space are given in
\cite{hl}, \cite{l1}, \cite{l2}, \cite{t}.  (An error in \cite{l1} is
corrected in \cite{t}.)  A regularity theorem for Einstein metrics has been
obtained by Skinner \cite{s}.

\section{Conformally Compact Einstein Metrics}

Let $X$ be the interior of a compact manifold with boundary $\Xb$ of
dimension $n+1$ as in the
introduction and let $g_+$ be a conformally compact metric on $X$.  Let $r$
be a sufficiently 
smooth defining function for $M=\pa X$ defined near $M$ and set 
$\gb = r^2 g_+$.  As discussed in the introduction, the conformal class 
$[\gb|_{TM}]$ is an invariant of $g_+$, independent of any choices.
The function $|dr|^2_{\gb} = \gb^{ij}r_ir_j$ extends
to $\Xb$ and its restriction to $M$ is independent of the choice
of $r$, so defines a second invariant of $g_+$.  The metric $g_+$ on $X$ is
complete and its sectional curvature is asymptotically constant at each
boundary point--conformally transforming the curvature tensor shows that 
\begin{equation}\label{curvten}
R_{ijkl} = -(|dr|^2_{\gb})(g_{ik}g_{jl}-g_{il}g_{jk}) +O(r^{-3}),
\end{equation}
where here the curvature tensor $R$ and metric $g$ both refer to $g_+$, and
our conventions are such that the above formula without the error term
defines a curvature tensor of constant curvature $-(|dr|^2_{\gb})$.  It
follows that the value of the invariant $|dr|^2_{\gb}$ at a boundary point
is the negative of the asymptotic sectional curvature of $g_+$ there.  

We will assume that $g_+$ satisfies the normalized Einstein condition
$\Ric(g_+)=-ng_+$.  Contracting in (\ref{curvten}) shows that in this case 
we have $|dr|^2_{\gb} = 1$ on $M$.  

In general, a choice of defining function $r$ determines a representative
metric $\gb|_{TM} = (r^2 g_+)|_{TM}$ for the conformal structure on $M$.
However, 
in the other direction, the conformal representative and this relation only
determine $r \mod O(r^2)$.  In the case when $|dr|^2_{\gb} = 1$ on $M$, in
particular when $g_+$ is Einstein, one can impose a second condition to
determine $r$ uniquely in a neighborhood of $M$.

\begin{lemma}\label{deffn}
A metric on $M$ in the conformal
infinity of $g_+$ determines a unique defining function $r$ in 
a neighborhood of $M$ such that $\gb|_{TM}$ is the prescribed
boundary metric and such that $|dr|^2_{\gb}=1$.
\end{lemma}
\begin{proof}
Given any choice of defining function $r_0$, 
let $\gb_0=r^2_0 g_+$ and set $r=r_0 e^{\om}$, so $\gb=e^{2\om}\gb_0$
and $dr=e^\om (dr_0+r_0d\om )$.  Thus
$$
|dr|^2_{\gb}=|dr_0 + r_0d\om |^2_{\gb_0}=
|dr_0|^2_{\gb_0}+2r_0(\nabla_{\gb_0}r_0)(\om )
+r_0^2|d\om |^2_{\gb_0},$$
so the condition $|dr|^2_{\gb}=1$ is equivalent to
\begin{equation}\label{pde}
2(\nabla_{\gb_0}r_0)(\om )+r_0|d\om |^2_{\gb_0}=
\frac{1-|dr_0|^2_{\gb_0}}{r_0}.
\end{equation}
This is a non-characteristic first order PDE for $\om$, so there is a
solution near $M$ with $\om |_{M}$ arbitrarily prescribed.
\end{proof}

A defining function determines for some $\epsilon>0$ an
identification of $M \times [0,\epsilon)$ with a neighborhood of
$M$ in $\Xb$: $(p,\lambda)\in M\times[0,\epsilon)$
corresponds to the point obtained by following the integral curve of
$\nabla_{\gb}r$ emanating from $p$ for $\lambda$ units of time.  For a
defining function of the type given in the lemma, with
$|dr|^2_{\gb}=1$, the $\lambda$-coordinate is just $r$, and
$\nabla_{\gb}r$ is orthogonal to the slices
$M\times\{\lambda\}$.  Hence, identifying $\lambda$ with $r$,
on $M\times[0,\epsilon)$ the metric $\gb$ takes the form
$\gb=g_{r}+dr^2$ for a 1-parameter family $g_{r}$ of metrics on
$M$, and 
\begin{equation}\label{form}
g_+=r^{-2}(g_{r}+dr^2).
\end{equation}

We explicitly identify a special defining function $r$ and normal form
(\ref{form}) for the hyperbolic metric 
$g_+ = 4(1-|x|^2)^{-2}\Sigma (dx^i)^2$ 
on the unit ball in $\R^{n+1}$.  Notice that in general the
condition $|dr|^2_{\gb}=1$ can be rewritten as 
$|d(\log \frac1r)|^2_{g_+}=1$, 
which is the eikonal equation for $\log \frac1r$ in the metric $g_+$.
The distance function $d(x) =$ (hyperbolic distance from $x$ to $0$)
satisfies the eikonal equation and also $d(x) \rightarrow \infty$ as
$|x| \rightarrow 1$, so we take $\log \frac1r = d(x)$, i.e. $r=e^{-d(x)}$.  
Now it is a basic fact of hyperbolic geometry that 
$d(x) = \log \frac{1+|x|}{1-|x|}$, so $r=\frac{1-|x|}{1+|x|}$ is a special
defining function for $\Hbb^{n+1}$ as in Lemma~\ref{deffn}.  Then
$\gb = r^2g_+=4(1+|x|)^{-4}\Sigma (dx^i)^2$, so the associated
representative for the conformal structure is 
$g_0 = \frac14(\mbox{usual metric on }\Sbb^n$).  Writing $\Sigma (dx^i)^2$
in polar coordinates and expressing everything in terms of $r$ gives
$g_+ = r^{-2}\left((1-r^2)^2g_0+ (dr)^2 \right)$, and therefore
\begin{equation}\label{hyp}
g_r = (1-r^2)^2g_0.
\end{equation}

We now impose the Einstein condition on a metric of the form (\ref{form}).
One can decompose the tensor $\Ric(g_+) + ng_+$ into components with respect
to the product structure $M\times (0,\epsilon)$.  
A straightforward calculation shows that the vanishing of
the component with both indices in $M$ is given by
\begin{equation}\label{ric}
rg_{ij}'' + (1-n)g_{ij}' - g^{kl}g_{kl}'g_{ij} - rg^{kl}g_{ik}'g_{jl}' 
+\frac r2g^{kl}g_{kl}'g_{ij}' - 2r\Ric_{ij}(g_r)=0,
\end{equation}
where $g_{ij}$ denotes the tensor $g_r$ on $M$, $'$ denotes $\pa_r$, and
$\Ric_{ij}(g_r)$ denotes the Ricci tensor of $g_r$ with $r$ fixed.  
As indicated in the introduction, we assume that 
$g_r$ is sufficiently regular that its asymptotics may be calculated 
from (\ref{ric}) (and the
equations for the other components of $\Ric(g_+) + ng_+$).  
Differentiating (\ref{ric}) $\nu - 1$ times with respect to $r$ and setting
$r=0$ gives
\begin{equation}\label{deriv}
(\nu-n)\pa_r^{\nu}g_{ij} - g^{kl}(\pa_r^{\nu}g_{kl})g_{ij} =
(\mbox{terms involving }\pa_r^{\mu}g_{ij}\mbox{ with } \mu<\nu).
\end{equation}
Beginning with the initial condition that $g_r$ is a given representative
metric at $r=0$, we may use (\ref{deriv}) inductively to solve for the
expansion of $g_r$.  So long as $\nu<n$, $\pa_r^{\nu}g|_{r=0}$ 
is uniquely determined at each step, 
and since the left-hand side of (\ref{ric}) respects parity in $r$, we have
$\pa_r^{\nu}g|_{r=0}=0$ for $\nu$ odd.  However this breaks
down for $\nu = n$.  In that case, if $n$ is odd, it follows from parity
considerations that the right-hand side of (\ref{deriv}) vanishes at $r=0$, 
so $g^{kl}\pa_r^ng_{kl}=0$ but the trace-free part of $\pa_r^ng_{kl}$ 
may be chosen arbitrarily.  If $n$ is even, then the 
right-hand side of (\ref{deriv}) might have non-vanishing
trace-free part, forcing the inclusion of a $r^n \log r$ term in the
expansion for $g_r$ with a trace-free coefficient.  The trace of
the $r^n$ coefficient is determined but not its trace-free part.
It can be shown that the remaining components of $\Ric(g_+)+ng_+$
give no further information to this order.

Summarizing, we see that for $n$ odd, the expansion of $g_r$ is of the form
\begin{equation}\label{oddasym}
g_r= g^{(0)} + g^{(2)} r^2 + (\mbox{even powers}) + g^{(n-1)} r^{n-1} +
g^{(n)} r^n + \ldots,
\end{equation}
where the $g^{(j)}$ are tensors on $M$, and $g^{(n)}$ is trace-free with
respect to a metric in the conformal class on $M$.  For $j$ even
and $0\leq j \leq n-1$, the tensor $g^{(j)}$ is locally formally determined
by the conformal representative, but $g^{(n)}$ is formally
undetermined, subject to the trace-free condition.  For $n$ even the
analogous expansion is
\begin{equation}\label{evenasym}
g_r= g^{(0)} + g^{(2)} r^2 + (\mbox{even powers}) + hr^n \log r +
g^{(n)} r^n + \ldots,
\end{equation}
where now the $g^{(j)}$ are locally determined for $j$ even and $0\leq j
\leq n-2$, $h$ is locally determined and trace-free, the trace of
$g^{(n)}$ is locally determined, but the trace-free part of $g^{(n)}$ is
formally undetermined. 

Of course, the determined coefficients in these expansions may be
calculated by 
carrying out the indicated differentiations above and keeping track of the
lower order terms at each stage.  For example, for $n=2$ one finds that 
$h=0$ and 
\begin{equation}\label{n=2}
g^{ij}g^{(2)}_{ij} = -\frac 12 R,
\end{equation}
while for $n\geq 3$ one has $g^{(2)}_{ij}=-P_{ij}$, where
\begin{equation}\label{P}
(n-2)P_{ij} = R_{ij} - \frac{R}{2(n-1)}g_{ij},
\end{equation}
and $R_{ij}$ and $R$ denote 
the Ricci tensor and scalar curvature of the chosen representative $g_{ij}$
of the conformal infinity.

In order to establish conformal invariance of the renormalized volume
invariants, we will later need to use the following Lemma.
\begin{lemma}\label{rel}
Let $r$ and $\rh$ be special defining functions as in Lemma~\ref{deffn}
associated to two different conformal representatives.  Then 
\begin{equation}\label{relation}
\rh = r e^{\om}
\end{equation} 
for a function $\om$ on $M\times[0,\ep)$ whose expansion at
$r=0$ consists only of even powers of $r$ up through and including the
$r^{n+1}$ term.
\end{lemma}
\begin{proof}
We have $\rh=e^{\om}r$ where $\om$ is determined by (\ref{pde}), which in
this case becomes 
\begin{equation}\label{newpde}
2\om_r +r(\om_r^2+|d_M\om|_{g_r}^2)=0.
\end{equation}  
The Taylor expansion of $\om$ is determined inductively by differentiating
this equation at $r=0$.  Clearly $\om_{r} = 0$ at $r=0$.  Consider
the determination of $\partial_r^{k+1}\om$ resulting from
differentiating (\ref{newpde}) an even number $k$ times and
setting $r=0$.  The term $\om_r^2$ gets differentiated $k-1$ times, 
so one of the two factors ends up differentiated an odd number of
times, so by induction vanishes at $r=0$.  Now
$|d_M\om|_{g_r}^2=g_{r}^{ij}\om_i\om_j$, so the $k-1$ differentiations
must be split between the three factors, so one of the factors 
must receive an odd
number of differentiations.  When an odd number of derivatives hits a $\om_i$,
the result again vanishes by induction.  But by (\ref{oddasym}) and
(\ref{evenasym}), so long as $k-1<n$, the odd derivatives of $g_r$ vanish
at $r=0$.  
\end{proof}

\section{Volume Renormalization}
Let $g_+$ be a conformally compact Einstein metric on $X$.  
As discussed above, a representative metric $g$ on $M$ for the 
conformal infinity of $g_+$ determines a special defining function $r$ for
$M$ and an identification
of a neighborhood of $M$ in $\Xb$ with $M \times [0,\ep)$.  In this
identification, 
$g_+$ takes the form (\ref{form}), where $g_0 = g$ is the chosen
representative metric.  Therefore the volume element $dv_{g_+}$ is given by 
\begin{equation}\label{volform}
dv_{g_+} = r^{-n-1}\left(\frac {\det g_r}{\det g}\right)^{1/2} dv_gdr.
\end{equation}
{From} (\ref{oddasym}) and (\ref{evenasym}) and the properties stated there 
for the coefficients in those expansions, it follows that
\begin{equation}\label{det}
\left(\frac {\det g_r}{\det g}\right)^{1/2}=
1+v^{(2)}r^2+(\mbox{even powers})+v^{(n)}r^{n}+ \ldots,
\end{equation}
where the $\ldots$ indicates terms vanishing to higher order.  All
indicated $v^{(j)}$ are locally determined functions on $M$ and 
$v^{(n)}=0$ if $n$ is odd.
 
Consider now the asymptotics of $\Vol_{g_+}(\{r>\ep\})$ as $\ep\rightarrow
0$.  Pick a small number $r_0$ and express
$\Vol(\{r>\ep\})= C + \int_{\{\ep<r<r_0\}}dv_{g_+}$.  
Integrating (\ref{volform}) using (\ref{det}) we obtain for $n$ odd
\begin{equation}\label{vodd}
\begin{array}l
\Vol(\{r>\ep\})= c_0 \ep^{-n} + c_2\ep^{-n+2} +
(\mbox{odd powers}) + c_{n-1} \ep^{-1} \\[3pt]
\phantom{\Vol(\{r>\ep\})=} + V + o(1)
\end{array}
\end{equation}
and for $n$ even
\begin{equation}\label{veven}
\begin{array}l
\Vol(\{r>\ep\})=c_0 \ep^{-n} + c_2\ep^{-n+2} +
(\mbox{even powers})+ c_{n-2} \ep^{-2} \\[3pt]
\phantom{\Vol(\{r>\ep\})=}  + L\log {\frac{1}{\ep}} +V + o(1). 
\end{array}
\end{equation}
The coefficients $c_i$ and $L$ are integrals over $M$ of local curvature
expressions of the metric $g$.  For example, 
$c_0 = \frac 1n \Vol_g(M)$.
Also,
\begin{equation}\label{Lform}
L = \int_M v^{(n)}\, dv_g.
\end{equation}

The renormalized volume is the constant term $V$ in the expansion for 
$\Vol(\{r>\ep\})$, which a-priori depends on the choice $g$ of
representative conformal metric on $M$.  

\begin{theorem}\label{invariance}
If $n$ is odd, then $V$ is independent of the choice of $g$.\\
If $n$ is even, then $L$ is independent of the choice of $g$.
\end{theorem}
\begin{proof}
The special defining functions $r$ and $\rh$ associated to
representative metrics $g$ and $\gh$ are related as in Lemma \ref{rel}. 
We can solve (\ref{relation}) for $r$ to give $r=\rh b (x,\rh)$, where the
expansion of $b$ also has only even powers of $\rh$ up through the
$\rh^{n+1}$ term.  It is important to note that in this relation, the $x$
still refers to the identification associated with $r$.

Set $\eh(x,\ep) = \ep b(x,\ep)$.  Then $\rh>\ep$ is equivalent to 
$r>\eh(x,\ep)$, so
$$
\Vol(\{r>\ep\})-\Vol(\{\rh>\ep\})= 
\int_{M}\int_{\ep}^{\eh}dv_{g_+}
$$
\begin{equation}\label{diff}
=\int_M\int_{\ep}^{\eh}
\sum_{\stackrel{0\leq j\leq n}{j\, {\rm even}
}} v^{(j)}(x) r^{-n-1+j}dr dv_g + o(1),
\end{equation}
where we have used (\ref{volform}), (\ref{det}).
For $n$ odd this is 
$$
\sum_{\stackrel{0\leq j\leq n-1}{j\, {\rm even}}}\ep^{-n+j}
\int_M\frac{v^{(j)}(x)}{-n+j}\left ( b(x,\ep)^{-n+j}-1 \right )dv_g +
o(1). 
$$
Since $b(x,\epsilon)$ is even through terms of order $n+1$ in $\ep$,
it follows that this expression has no constant term as
$\epsilon \rightarrow 0$.  Similarly, when $n$ is even, the $r^{-1}$ term 
in (\ref{diff}) contributes $\log b(x,\ep)$, so there is no
$\log {\frac{1}{\ep}}$ term as $\epsilon \rightarrow 0$. 
\end{proof}

According to Theorem \ref{invariance}, for $n$ odd the renormalized volume
$V$ is an absolute invariant of the conformally compact Einstein metric
$g_+$.  But this is not so if $n$ is even.  Suppose $g$ and $\gh =
e^{2\up}g$ are two metrics in the conformal infinity of $g_+$, where $\up
\in C^{\infty}(M)$.  The difference $V_g - V_{\gh}$
is the constant term in the expansion of (\ref{diff}).  By the local
determination of the $v^{(j)}$ and of the expansion of $b(x,\ep)$, we see
that this anomaly takes the form
$$
V_{\gh}-V_g=\int_M \mathcal P_g(\up)dv_g,
$$
where $\mathcal P_g$ is a polynomial nonlinear differential operator whose
coefficients are polynomial expressions in $g$, its inverse, and its
derivatives.  Moreover, it is easy to see that 
the linear part in $\up$ of $\mathcal P_g(\up)$ is just
$v^{(n)}\up$.  Since this linear part measures the infinitesimal change
under conformal rescalings,  $V_{\gh}-V_g$ is determined
by knowledge of  $v^{(n)}$ for general $g$.  In summary, for $n$ even,
the fundamental object is the function $v^{(n)}$---its integral over $M$ is
by (\ref{Lform}) the conformal invariant $L$, and 
multiplication by it gives the infinitesimal
anomaly, which determines the full anomaly.

It is straightforward to carry out the calculations indicated above to
identify $v^{(n)}$ and $\mathcal P_g$ in low dimensions.  For $n=2$ one
obtains  
$$v^{(2)} = -\frac 14 R, \qquad \mathcal P_g(\up)=-\frac14 (R\up +
\up_i\up^i),$$ 
so $L=-\pi \chi(M)$, where $\chi(M)$ denotes the Euler characteristic of
$M$. 

For $n=4$ one obtains
$$v^{(4)}= \frac18 [(P_i{}^i)^2 - P_{ij}P^{ij}],$$
$$\mathcal P_g(\up) = v^{(4)}\up + \up_{ij}\up^i\up^j - P_{ij}\up^i\up^j - 
\frac14(\up_i\up^i)^2 +P_j{}^j\up_i\up^i.$$
The Gauss-Bonnet Theorem for $n=4$ reads 
$$32\pi^2 \chi(M) = \int_M [|W|^2 -8P_{ij}P^{ij} + 8(P_i{}^i)^2] 
dv_g,$$ where 
$$W_{ijkl}=R_{ijkl} - (P_{ik}g_{jl}
+P_{jl}g_{ik}-P_{il}g_{jk}-P_{jk}g_{il})$$ 
is the Weyl conformal curvature tensor.  Therefore
$$L= \frac {\pi^2}{2} \chi(M) - \frac {1}{64} 
\int_M |W|^2 dv_g.$$

For $n=6$ one obtains
$$v^{(6)}= \frac {1}{48}[-P^{ij}B_{ij} 
+3P_i{}^iP_{kl}P^{kl}-2P_{ij}P_k{}^iP^{jk} -(P_i{}^i)^3],$$
where
$$B_{ij}=P_{ij,k}{}^k -P_{ik,j}{}^k - P^{kl}W_{kijl}.$$
Again there is an explicit realization of $L=\int_M  v^{(6)}dv_g$ 
as a linear
combination of the Euler characteristic and the integral of a local
conformal invariant.  Define 
$$C_{ijk}= P_{ij,k}-P_{ik,j}$$
and set
$$V_{ijklm}=W_{ijkl,m}+g_{im}C_{jkl}-g_{jm}C_{ikl}+g_{km}C_{lij}
-g_{lm}C_{kij}$$ and
$$U_{ijkl}=C_{jkl,i}-P_i{}^mW_{mjkl}.$$
Then $$I=|V|^2 - 16 W_{ijkl}U^{ijkl} +16 |C|^2$$ 
is a conformal invariant in general dimension $n\geq3$; it is the
norm-squared of the first covariant derivative of the 
curvature tensor of the ambient metric of \cite{fg}. 
One can calculate that for $n=6$,
$$L=-\frac{\pi^3}{6} \chi(M) + \frac {1}{2304}\int_M J dv_g,$$
where 
$$J = -3I 
+7W_{ijkl}W^{ij}{}_{pq}W^{klpq} +4 W_{ijkl}W^{ipkq}W^j{}_p{}^l{}_q.$$

For $\Hbb^{n+1}$, using (\ref{hyp}) it is possible to calculate the
invariants $V$ for $n$ odd and $L$ for $n$ even.  
{From} (\ref{hyp}) one obtains
$$\left(\frac {\det g_r}{\det g_0}\right)^{1/2} = (1-r^2)^n,$$ 
so recalling that $4g_0$ is the usual metric on $\Sbb^n$, it follows from 
(\ref{volform}) that
\begin{equation}\label{hypvol}
\Vol(\{r>\ep\}) = 2^{-n}\mbox{Area}(\Sbb^n)\int_{\ep}^1 r^{-n-1}
(1-r^2)^n dr.
\end{equation}

For $n$ odd, write
$$\int_{\epsilon}^1 r^{-n-1}(1-r^2)^ndr = -\frac{1}{n}
\int_{\epsilon}^1(1-r^2)^nd(r^{-n}) $$
$$= 
\frac{1}{n}\epsilon^{-n}(1-\epsilon^2)^n-
2\int_{\epsilon}^1r^{-n+1}(1-r^2)^{n-1}dr.$$
The boundary term has no constant term in $\epsilon$, 
so upon applying the same procedure
inductively it follows that $\int_{\epsilon}^1r^{-n-1}(1-r^2)^ndr$
has constant term 
$$\frac{(-2)^{\frac{n+1}{2}}n(n-1)\ldots \left(\frac{n+1}{2}\right)}
{n(n-2)\ldots 1}\int_0^1(1-r^2)^{\frac{n-1}{2}}dr.$$
Collecting the constants, one finds 
$$V=(-1)^{\frac{n+1}{2}}
\frac{\pi^{\frac{n+2}{2}}}{\Gamma(\frac{n+2}{2})}.$$ 

For $n=2m$ even, expand $(1-r^2)^n$ using the binomial theorem; it follows
that the $\log \frac {1}{\ep}$ coefficient in the expansion of 
$\int_{\ep}^1 r^{-n-1}(1-r^2)^n dr$
 is $(-1)^{m}
\left(\begin{array}{c} n \\ m \end{array} \right)$.  Substituting into 
(\ref{hypvol}) and simplifying gives 
$$L=(-1)^{m}\frac {2\pi^{m}}{m!}.$$

A more familiar setting for conformal anomalies is in the study of
functional determinants of conformally invariant differential operators.
The invariance properties of $V$ are reminiscent of those for the
functional determinant of the conformal Laplacian, which is conformally
invariant in odd dimensions but which has an anomaly in even
dimensions (\cite{pr}).  We remark that the AdS/CFT
correspondence predicts that the volume anomaly for $n=4$ is a particular
linear combination of functional determinant anomalies on scalars, spinors,
and 1-forms; this prediction was confirmed in \cite{hs1}. 
The properties of the invariant $L$ are, on the other hand, similar to
those for the constant term in the
expansion of the integrated heat kernel for the conformal Laplacian,
which vanishes in odd dimensions but in even
dimensions is a conformal invariant obtained by integrating a local
expression in curvature (\cite{bo}, \cite{pr}).  

\section{Area Renormalization}
Let $(X^{n+1},g_+)$ be a conformally compact Einstein manifold with
conformal infinity $(M,[g])$ as above.  In this section we describe the
renormalization of the area 
of minimal submanifolds $Y \subset X$ of dimension $k+1$, 
$0 \leq k \leq n-1$, which extend regularly to
$\Xb$.  Set $N=\Yb \cap M$.  We assume that $N$ is a smooth submanifold of 
$M$.  We will outline the arguments and refer to \cite{gw} for details.

First one must study the asymptotics of $Y$ near $M$.
Locally near a point of $N$, coordinates $(x^{\al},u^{\al'})$ 
for $M$ may be chosen, where $1\leq \al \leq k$ and
$1\leq \al' \leq n-k$, so that $N=\{u=0\}$ and so that $\partial_{x^{\al}}
\perp \partial_{u^{\al'}}$ on $N$ with respect to 
a metric in the conformal infinity of $g_+$.  Choose a representative
metric $g$ for the conformal infinity and recall that this choice
determines by Lemma \ref{deffn} a defining function $r$ for $M$ and an
identification of a neighborhood of $M$ in $\Xb$ with $M\times [0,\ep)$.  
This identification determines an extension of the $x^{\al}$ and $u^{\al'}$ 
into $X$, and together with $r$ these form a local 
coordinate system on $\Xb$.  We
consider submanifolds $Y$ which in such coordinates may be written as a 
graph $\{u=u(x,r)\}$.  One can calculate the minimal surface equation
for $Y$ explicitly as a system of differential equations for the unknowns 
$u^{\al'}(x,r)$.  These equations together with the boundary condition
$u(x,0)=0$ are used to study the asymptotics of $u(x,r)$ at $r=0$.  One
finds (see \cite{gw}) that for $k$ odd
\begin{equation}\label{uodd}
u= u^{(2)} r^2 + (\mbox{even powers}) + u^{(k+1)} r^{k+1} +
u^{(k+2)} r^{k+2} + \ldots,
\end{equation}
and for $k$ even
\begin{equation}\label{ueven}
u= u^{(2)} r^2 + (\mbox{even powers}) + u^{(k)} r^{k} +
wr^{k+2} \log r + u^{(k+2)} r^{k+2} + \ldots,
\end{equation}
where the $u^{(j)}$ and $w$ are functions of $x$, all of which are locally 
determined except for $u^{(k+2)},$ and the $\ldots$ indicates terms
vanishing to higher order.  Observe in particular that the minimal
submanifold $Y$ is determined to order $k+2$ by $N = \partial Y$, that 
the expansion of $u$ is even in $r$ to order $k+2$, and that the
irregularity in the expansion occurs at order $k+2$.  The consequence
$\pa_r u=0$ at $r=0$ has the geometric interpretation 
that $Y$ intersects $M$ orthogonally, a fact
very familiar from the geometry of geodesics in hyperbolic space.  For
the case $k=0$ of geodesics it turns out that necessarily $w=0$, and the
local indeterminancy in this case of $u^{(2)}$ is a reflection of the
familiar fact 
that at the boundary a geodesic may have any asymptotic curvature measured
with respect to the smooth metric $\gb$.

Next one calculates the metric induced on $Y$ by the
conformally compact Einstein metric $g_+$.  The area form 
$da_Y$ of $Y$ takes the form
\begin{equation}\label{areaform}
da_Y=r^{-k-1}\left [ 1 + a^{(2)} r^2 + (\mbox{even powers}) + a^{(k)}
r^{k} +\ldots \right ] da_Ndr,
\end{equation}
where the $\ldots$ indicates terms vanishing to higher order and $da_N$
denotes the area form on $N$ with respect to the chosen conformal
representative $g$ on the boundary.  All indicated $a^{(j)}$ are locally
determined functions on $N$ and $a^{(k)} =0$ if $k$ is odd.  A key
observation in establishing (\ref{areaform}) is
that since the induced metric depends only on $u$ and its first
coordinate derivatives, the local indeterminacy and irregularities 
at order $k+2$ in $u$ and those at order $n$ in the metric $g_r$ given by 
(\ref{oddasym}), (\ref{evenasym})
do not enter into the asymptotics of the area form to the indicated order.
The evenness of $r^{k+1}da_Y$ then follows from that of $g_r$ and of $u$. 

Now we can consider the asymptotics of Area$_{g_+}(Y \cap \{r>\ep\})$ as 
$\ep \rightarrow 0.$  Pick a small
number $r_0$ and express
$\mbox{Area}(Y \cap \{r>\ep\})= C + \int_{Y \cap \{\ep<r<r_0\}}da_Y$.  
By (\ref{areaform}) we obtain for $k$ odd
$$
\begin{array}l
\mbox{Area}(Y \cap \{r>\ep\})= b_0 \ep^{-k} + b_2\ep^{-k+2} +
(\mbox{even powers}) + b_{k-1} \ep^{-1} \\[3pt]
\phantom{\mbox{Area}(Y \cap \{r>\ep\})=} + A + o(1)
\end{array}
$$
and for $k$ even
\begin{equation}\label{aeven}
\begin{array}l
\mbox{Area}(Y \cap \{r>\ep\})= b_0 \ep^{-k} + b_2\ep^{-k+2} +
(\mbox{even powers}) + b_{k-2} \ep^{-2} \\[3pt]
\phantom{\mbox{Area}(Y \cap \{r>\ep\})=} + K\log {\frac{1}{\ep}} +A + o(1).  
\end{array}
\end{equation}
Observe that 
\begin{equation}\label{Kform}
K = \int_N a^{(k)}\, da_N.
\end{equation}
The analogue of Theorem \ref{invariance} is the following, which is proved by
a similar argument.
\begin{theorem}\label{ainvariance}
If $k$ is odd, then $A$ is independent of the choice of $g$.\\
If $k$ is even, then $K$ is independent of the choice of $g$.
\end{theorem}
Therefore, for $k$ odd, a minimal submanifold of $X$ has a well-defined
invariant renormalized area $A$.  For $k$ even, the $\log
\frac {1}{\ep}$ coefficient $K$ is a conformal invariant of the submanifold
$N$ of $M$ given according to (\ref{Kform}) by the integral of an
expression determined locally by the geometry of $N\subset M$ with respect
to the metric $g$.

Analogously to the volume case, there is a conformal anomaly 
for $A$ when $k$ is even.  If
$\gh =e^{2\up}g$ is a conformally related metric, then the local
determination of the coefficients $a^{(j)}$ in (\ref{areaform}) and of the
defining function $\rh$ as in Lemma \ref{deffn} implies that 
$$
A_{\gh}-A_g=\int_N \mathcal Q_N(\up)da_N
$$
for a differential expression $\mathcal Q_N$ determined locally by the
geometry of $N\subset M$.  One interesting difference from the volume
anomaly is that the linearization of $\mathcal Q_N(\up)$ need not be just
$a^{(k)}\up$--it can in general involve derivatives of $\up$ as well.
However it is clear from rescaling in (\ref{aeven}) that 
$\mathcal Q_N(\up) = a^{(k)}\up$ for $\up$ constant.

The invariant $K$ and the anomaly for the
lowest dimensional cases $k=0,2$ are calculated in \cite{gw}.  
For $k=0$, $Y$ is a union of geodesics
in $X$ and $N$ consists of finitely many points.  Of course a point has
no geometry and the conclusions are rather trivial; 
one finds that $K$ is the number of
boundary points, $\mathcal Q$ evaluates $\up$ at a boundary point,
and the anomaly is given by $A_{\gh}-A_g = \sum_{p\in N} \up(p)$.
To describe the $k=2$ results recall that the second fundamental form of 
$N\subset M$ with respect to the metric $g$ is the symmetric form 
$B^{\ga'}_{\al\be}$ on $TN$ with values in $TN^{\perp}$ defined by
$B(X,Y) = (\nabla_X Y)^{\perp}$ for
vectors $X,Y \in TN$; here $\nabla$ denotes the Levi-Civita covariant
derivative of $g_{ij}$ and $\perp$ the component in $TN^{\perp}$.
On $N$, the metric $g_{ij}$ decomposes into two pieces 
$g_{\al\be}$ and $g_{\al'\be'}$.
The mean curvature vector of $N$ is $H^{\ga'} = 
g^{\al\be} B^{\ga'}_{\al\be}$.
The tensor $P$ given by (\ref{P}) also decomposes into pieces with respect
to the decomposition $TM = TN \oplus (TN)^{\perp}$; we denote by 
$P_{\al \be}$ its component with both indices in $TN$ 
({\em not} the corresponding tensor for
the induced metric $g_{\al \be}$).  Then for $k=2$ one finds
\begin{equation}\label{K}
K=-\frac 18 \int_N(|H|^2 + 4g^{\al \be}P_{\al \be})da_N
\end{equation}
and
$$
\mathcal Q_N(\up)= -\frac 18 (|H|^2 + 4g^{\al \be}P_{\al \be})\up 
+\frac 14(H^{\ga'}\up_{\ga'} -\up_i\up^i) .
$$

The quantity defined by $(\ref{K})$ is therefore a conformal invariant of a
surface $N$ in a conformal manifold $M$.  For conformally flat space this
reduces to a multiple of 
the Willmore functional (for which, see, e.g., \cite{b}).  Other 
generalizations of the Willmore functional to curved conformal spaces are
given in \cite{c} and \cite{w}.

A different conformal anomaly associated to a surface in a conformal
6-manifold is discussed in \cite{hs2}.

\end{document}